%
%
%
%
%
\font\tenbf=cmbx10

\font\eightrm=cmr8
\font\eightit=cmti8
\font\germ=eufm10

\def\s{\hbox{\germ S}}
\def\sectiontitle#1\par{\vskip0pt plus.1\vsize\penalty-250
\vskip0pt plus-.1\vsize\bigskip\vskip\parskip
\message{#1}\leftline{\tenbf#1}\nobreak\vglue 5pt}
\def\ds{\displaystyle}

\def\eno{\eqalignno}
\def\ld{\lambda}

\def\frac#1#2{{#1\over#2}}

\magnification=\magstep1
\vsize=8.7truein
\parindent=15pt
\nopagenumbers
\baselineskip=10pt
\line{ 
\hfil}
\vglue 5pc
\baselineskip=13pt
\headline{\ifnum\pageno=1\hfil\else%
{\ifodd\pageno\rightheadline \else \leftheadline\fi}\fi}
\def\rightheadline{\hfil\eightit 
Skew divided difference operators and Schubert polynomials
\quad\eightrm\folio}
\def\leftheadline{\eightrm\folio\quad 
\eightit 
Anatol N. Kirillov 
\hfil}
\voffset=2\baselineskip
\centerline{\tenbf 
SKEW  \hskip 0.1cm DIVIDED \hskip 0.1cm DIFFERENCE \hskip 0.1cm
OPERATORS}
\vskip 0.1cm  
\centerline{\tenbf AND\hskip 0.1cm SCHUBERT \hskip 0.1cm POLYNOMIALS}
\vglue 12pt
\centerline{\eightrm 
ANATOL N. KIRILLOV
}
\vglue 6pt
\baselineskip=12pt
\centerline{\eightit
CRM, University of Montreal
}
\baselineskip=10pt
\centerline{\eightit 
C.P. 6128, Succursale A, Montreal (Quebec) H3C 3J7, Canada 
}
\baselineskip=12pt
\centerline{\eightit
and }
\baselineskip=12pt
\centerline{\eightit 
Steklov Mathematical Institute,
}
\baselineskip=10pt
\centerline{\eightit 
Fontanka 27, St.Petersburg, 191011, Russia
}
\vglue 18pt
\centerline{  }
\centerline{\eightrm ABSTRACT}
\vglue 4pt
{\rightskip=1.5pc
\leftskip=1.5pc
\eightrm\parindent=1pc 
We study an action of the skew divided difference operators on the 
Schubert polynomials and give an explicit formula 
for structural constants for the Schubert polynomials in 
terms of certain weighted paths in the Bruhat order on the symmetric group.
We also prove that, under certain assumptions, the skew divided 
difference operators transform Schubert polynomials into polynomials 
with positive integer coefficients.
} 
\vglue12pt
\baselineskip=13pt
\overfullrule=1pt
\def\qed{\hfill$\vrule height 2.5mm width 2.5mm depth 0mm$}
\vskip 0.5cm
{\bf \S 0. Introduction.}
\vskip 0.3cm

In this paper we study skew divided difference operators with 
applications to the Littlewood--Richardson problem in the Schubert 
calculus. The Littlewood--Richardson problem in the Schubert calculus is 
a problem of finding a combinatorial rule for computation of the Schubert
polynomials structural constants $c^w_{uv}$, $u,v,w\in S_n$. 
The coefficients $c^w_{uv}$ are 
defined via decomposition of the product of two Schubert polynomials 
$\s_u$ and $\s_v$, $u,v\in S_n$,
$$\s_v\s_u\equiv\sum_{w\in S_n}c_{uv}^w\s_w~({\rm mod}~I_n), \eqno (0.1)
$$
where $I_n$ is the ideal in the ring of polynomials ${\bf Z}[x_1,\ldots 
,x_n]$ generated by the elementary symmetric functions without constant 
term. By now such rule is known only in the case when $u,v,w$ are the 
Grassmannian permutations (see, e.g., [M1], p.13, and [M2], Chapter~I, 
\S 9); this is the famous Littlewood--Richardson rule for Schur 
functions.

The skew divided difference operators were introduced by I.~Macdonald in 
[M1]. The simplest way to define the skew divided difference operators is 
based on the Leibnitz rule for the divided difference operators $\partial_w$, 
$w\in S_n$, namely,
$$\partial_w(fg)=\sum_{w\succeq 
v}\left(\partial_{w/v}f\right)\partial_vg, \eqno (0.2)
$$
The symbol $w\succeq v$ for $w,v\in S_n$, 
here and after, means that $w$ dominates $v$ with 
respect to the Bruhat order on the symmetric group $S_n$ (see, e.g., 
[M1], p.6). Formula (0.2) is reduced to the classical Leibnitz rule 
in the case when $w=(i,i+1)$ is a simple transposition:
$$\partial_{i}(fg)=\partial_{i}(f)g+s_{i}(f)\partial_{i}g.
$$
One of the main applications of the skew divided difference operators 
is an elementary and transparent algebraic proof of 
the Monk formula for Schubert polynomials (see [M1], (4.15)).

Our interest to the skew divided difference operators is based on their 
connection with the structural constants for Schubert polynomials. 
More precisely, if $w,v\in S_n$, and $w\succeq v$, then
$$\partial_{w/v}\left(\s_u\right)|_{x=0}=c^w_{uv}. \eqno (0.3)
$$
Polynomial 
$\partial_{w/v}(\s_u)$ is a homogeneous polynomial in $x_1,\ldots ,x_n$ 
of degree $l(u)+l(v)-l(w)$ with integer coefficients. We make a 
conjecture that in fact
$$\partial_{w/v}(\s_u)\in{\bf N}[x_1,\ldots ,x_n], \eqno (0.4)
$$
i.e. the polynomial $\partial_{w/v}(\s_u)$ has nonnegative integer 
coefficients. 
In the case $l(u)+l(v)=l(w)$, this conjecture follows from 
geometric interpretation of the structural constants $c^w_{uv}$ as the 
intersection numbers for Schubert cycles. For general $u,v,w\in S_n$ the 
conjecture is still open. In Section~7 we prove conjecture (0.4) in the
simplest, nontrivial case (see Theorem~1) 
when $w$ and $v$ are connected by an edge in the Bruhat order on the 
symmetric group $S_n$. In other words, if $w=vt_{ij}$, where $t_{ij}$ is 
the transposition that interchanges $i$ and $j$, and $l(w)=l(vt_{ij})+1$.
It is well--known ([M1], p.30) that in this case the skew divided 
difference operator $\partial_{w/v}$ coincides with operator $\partial_{ij}$, 
i.e. $\partial_{w/v}=\partial_{ij}$.
Our proof employs the generating function for Schubert polynomials 
("Schubert expression", [FS], [FK]) in the nilCoxeter algebra. 

In Section~8 we consider another 
application of the skew divided difference operators, namely, we give an 
explicit 
(but still not combinatorial) formula for structural constants for Schubert 
polynomials in terms of weighted paths in the Bruhat order with weights 
taken from the nilCoxeter algebra (see Theorem~2).

\vskip 0.2cm

{\bf Acknowledgments.}
This note is based on the lectures "Schubert polynomials" 
delivered in the Spring 1995 at the University of Minneapolis and in the 
Spring 1996 at the University of Tokyo. I would like to thank my colleagues 
from these universities 
for hospitality and support. My special thanks to Victor Reiner for 
numerous fruitful suggestions and discussions about Schubert polynomials.
\vskip 0.5cm

{\bf \S 1. Skew Schur functions.}
\vskip 0.3cm

In this Section we review the definition and basic properties of the skew 
Schur functions. For more details and proofs, see [M2], Chapter I, \S 5.
The main goal of this Section is to arise a problem of constructing skew 
Schubert polynomials with properties "similar" to the those for skew 
Schur functions (see properties (1.2)--(1.5) below).

Let $X_n=(x_1,\ldots ,x_n)$ be a set of independent variables, and $\ld , 
\mu$ be partitions, $\mu\subset\ld$, $l(\ld )\le n$.
\vskip 0.3cm

{\bf Definition 1.} {\it The skew Schur function $s_{\ld /\mu}(X_n)$ 
corresponding to the skew shape $\ld -\mu$ is define to be
$$s_{\ld /\mu}(X_n)=\det\left(h_{\ld_i-\mu_j-i+j}\right)_{1\le i,j\le n}, 
\eqno (1.1)
$$
where $h_k:=h_k(X_n)$ is the complete homogeneous symmetric function of 
degree $k$ in the variables $X_n=(x_1,\ldots ,x_n)$.}
\vskip 0.2cm

Below we list the basic properties of skew Schur functions:

a) Combinatorial formula:
$$s_{\ld /\mu}(X_n)=\sum_Tx^{w(T)}, \eqno (1.2)
$$
where summation is taken over all semistandard tableaux $T$ of the 
shape $\ld -\mu$ with 
entries not exceeding $n$; here $w(T)$ is the weight of tableau $T$
(see, e.g., [M2], p.5), and 
$x^{w(T)}:=x_1^{w_1}x_2^{w_2}\cdots x_n^{w_n}$.

b) Connection with structural constants for Schur functions:
$$s_{\ld /\mu}=\sum_{\nu ,~l(\nu )\le n}c_{\mu\nu}^{\ld}s_{\nu}, \eqno (1.3)
$$
where the coefficients $c_{\mu\nu}^{\ld}$ (the structural constants, or 
the Littlewood--Richardson numbers) are defined through the decomposition
$$s_{\mu}s_{\nu}=\sum_{\ld}c_{\mu\nu}^{\ld}s_{\ld}. \eqno (1.4)
$$

c) Littlewood--Richardson rule:
$$s_{\ld /\mu}=\sum_{\nu ,l(\nu )\le n} |{\rm Tab}^0(\ld -\mu ,\nu 
)|s_{\nu}, \eqno (1.5)
$$
where ${\rm Tab}^0(\ld -\mu ,\nu )$ is the set of all semistandard 
tableaux $T$ of shape $\ld -\mu$ and weight $\nu$ such that the 
corresponding reading word w$(T)$ of tableaux $T$ (see, e.g., [M2], 
Chapter~I, \S 9) is a lattice word (ibid). Thus,
$${\rm Mult}_{V_{\ld}}(V_{\nu}\otimes V_{\mu})=c_{\nu\mu}^{\ld}=
|{\rm Tab}^0(\ld -\mu ,\nu )|. \eqno (1.6)
$$
\vskip 0.5cm

{\bf \S 2. Divided difference operators.}
\vskip 0.3cm

{\bf Definition 2.} {\it Let $f$ be a function of $x$ and $y$ (and possibly 
other variables), the divided difference operator $\partial_{xy}$ is 
defined to be}
$$\partial_{xy}f={f(x,y)-f(y,x)\over x-y}.\eqno (2.1)
$$
\vskip 0.2cm

The operator $\partial_{xy}$ takes polynomials to polynomials and has 
degree $-1$. On a product $fg$, $\partial_{xy}$ acts according to the 
Leibnitz rule
$$\partial_{xy}(fg)=(\partial_{xy}f)g+(s_{xy}f)(\partial_{xy}g), \eqno 
(2.2)
$$
where $s_{xy}$ interchanges $x$ and $y$.

It is easy to check the following properties of divided difference 
operators $\partial_{xy}$:

a) ~~$\partial_{xy}s_{xy}=-\partial_{xy}$, 
$s_{xy}\partial_{xy}=\partial_{xy}$,

b) ~~$\partial_{xy}^2=0$,

c) ~~$\partial_{xy}\partial_{yz}\partial_{xy}=
\partial_{yz}\partial_{xy}\partial_{yz}$,

d) ~~$\partial_{xy}\partial_{yz}=\partial_{xz}\partial_{xy}+
\partial_{yz}\partial_{xz}$.

The next step is to define a family of divided difference operators 
$\partial_i$, $1\le i\le n-1$, action on the ring of polynomials with $n$
variables.

Let $x_1,x_2,\ldots ,x_n$ be independent variables, and let
$$P_n={\bf Z}[x_1,\ldots ,x_n].
$$
For each $i$, $1\le i\le n-1$, let
$$\partial_i=\partial_{x_i,x_{i+1}},
$$
be divided difference operator corresponding to the simple transposition 
$s_i=(i,i+1)$ which interchanges $x_i$ and $x_{i+1}$.

Each $\partial_i$ is a linear operator on $P_n$ of degree $-1$. The divided 
difference operators $\partial_i$, $1\le i\le n-1$, satisfy the following 
relations
$$\eno{
i)~~~~ &\partial_i^2=0,~~{\rm if}~~1\le i\le n-1,\cr
ii)~~~&\partial_i\partial_j=\partial_j\partial_i ,~~{\rm if}~~1\le i,j\le 
n-1,~~{\rm and}~~|i-j|>1,\cr
iii)~~&\partial_i\partial_{i+1}\partial_i=
\partial_{i+1}\partial_i\partial_{i+1},~~{\rm if}~~1\le i\le n-2.}
$$

Let $w\in S_n$ be a permutation; then $w$ can be written as a product of 
simple transpositions $s_i=(i,i+1)$, $1\le i\le n-1$, namely,
$$w=s_{i_1}\cdots s_{i_p}.
$$
Such a representation is called reduced if $p=l(w)$, where $l(w)$ is the 
length of $w$. For each $w\in S_n$, let $R(w)$ denote the set of all 
reduced decompositions of $w$, i.e. the set of all sequences $(i_1,\ldots 
,i_p)$ of length $p=l(w)$ such that $w=s_{i_1}\cdots s_{i_p}$.

For any sequence ${\bf a}=(a_1,\ldots ,a_p)$ of positive integers, let us 
define $\partial_{\bf a}=\partial_{a_1}\cdots\partial_{a_p}$.
\vskip 0.3cm

{\bf Proposition 1} ([M1], Chapter II). {\it

i) If a sequence ${\bf a}=(a_1,\ldots ,a_p)$ is not reduced, then 
$\partial_{\bf a}=0$.

ii) If ${\bf a,b}\in R(w)$ then $\partial_{\bf a}=\partial_{\bf b}$.}
\vskip 0.2cm

From Proposition~1, $ii)$ follows that one can define 
$\partial_w=\partial_{\bf a}$ unambiguously, where ${\bf a}$ is any 
reduced decomposition for $w$.
\vskip 0.5cm

{\bf \S 3. Schubert polynomials.}
\vskip 0.3cm

In this section we recall the definition and basic properties of the 
Schubert polynomials introduced by A.~Lascoux and M.-P.~Sch\"utzenberger. 
Further details and proofs can be found in [M1].

Let $\delta =\delta_n=(n-1,n-2,\ldots ,1,0)$, so that
$$x^{\delta}=x^{n-1}_1x^{n-2}_2\cdots x_{n-1}.
$$
\vskip 0.3cm

{\bf Definition 3} (Lascoux-Sch\"utzenberger, [LS]). {\it For each 
permutation $w\in S_n$ the Schubert polynomial $\s_w$ is defined to be
$$\s_w=\partial_{w^{-1}w_0}(x^{\delta}), \eqno (3.1)
$$
where $w_0$ is the longest element of $S_n$.}
\vskip 0.2cm

{\bf Proposition 2} {[LS], [M1]) {\it

i) For each permutation $w\in S_n$, $\s_w$ is a polynomial in $x_1,\ldots 
,x_{n-1}$ of degree $l(w)$ with positive integer coefficients.

ii) Let $v,w\in S_n$. Then}
$$\partial_v\s_w=\cases{\s_{wv^{-1}}, & if~ $l(wv^{-1})=l(w)-l(v)$,\cr
0, & otherwise.}
$$

{\it iii) Schubert polynomials $\s_w$, $w\in S_n$, form a ${\bf Z}$--linear 
basis in the space ${\cal F}_n$, where
$${\cal F}_n=\left\{ f\in 
P_n~|~f=\sum_{\alpha\subset\delta}c_{\alpha}x^{\alpha}\right\} .
$$

iv) Schubert polynomials $\s_w$, $w\in S_n$, form an orthogonal basis 
with respect to the pairing $\langle ~,~\rangle_0$:}
$$\langle\s_w,\s_u\rangle_0
=\cases{1, & if $u=w_0w$,\cr 0, & otherwise,}
$$
{\it where $\langle f,g\rangle_0=\eta 
(\partial_{w_0}(fg)):=\partial_{w_0}(fg)|_{x=0}$.

v) (Stability) Let $m>n$ and let $i:S_n\hookrightarrow S_m$ be the natural 
embedding. Then}
$$\s_w=\s_{i(w)}.
$$
\vskip 0.5cm

{\bf \S 4. Skew divided difference operators.}
\vskip 0.3cm

The skew divided difference operators $\partial_{w/v}$, $w,v\in S_n$, were 
introduced by I.~Macdonald, [M1], Chapter II.

Let $w,v\in S_n$, and $w\succeq v$ with respect to the Bruhat order 
$\succeq$ on the symmetric group
$S_n$. In other words, if ${\bf a}=(a_1,\ldots ,a_p)$ is a reduced 
decomposition for $w$ then there exists a subsequence ${\bf b}\subset{\bf 
a}$ such that ${\bf b}$ is a reduced decomposition for $v$ (for more 
details, see, e.g., [M1], (1.17)).
\vskip 0.3cm

{\bf Definition 4} (Macdonald, [M1]). {\it Let $v,w\in S_n$, and $w\succeq 
v$ with respect to the Bruhat order, and ${\bf a}=(a_1,\ldots a_p)\in 
R(w)$. The skew divided difference operator $\partial_{w/v}$ is defined 
to be
$$\partial_{w/v}=v^{-1}\sum_{{\bf b}\subset{\bf a},~{\bf b}\in R(v)}
\phi ({\bf a},{\bf b}), \eqno (4.1)
$$
where}
$$\eno{
&\phi ({\bf a},{\bf b})=\prod_{i=1}^p\phi_i({\bf a},{\bf b}),\cr
&\phi_i({\bf a},{\bf b})=\cases{s_{a_i}, & if $a_i\in{\bf b}$, \cr 
\partial_{a_i}, & if $a_i\not\in{\bf b}$ .}}
$$
\vskip 0.2cm

One can show (see, e.g., [M1], p.29) that Definition~4.1 is independent 
of the reduced decomposition ${\bf a}\in R(w)$.

Below we list the basic properties of the skew divided difference 
operators $\partial_{w/v}$. For more details and proofs, see, e.g., [M1].
The statement iv) of Proposition~3 below seems to be new.
\vskip 0.3cm

{\bf Proposition 3.} {\it

i) Let $f,g\in P_n$, $w\in S_n$, then
$$\partial_w(fg)=\sum_{w\succeq v}(\partial_{w/v}f)\partial_vg. \eqno (4.2)
$$

More generally,

ii) Let $f,g\in P_n$, $u,w\in S_n$, and $w\succeq u$ with respect to the the 
Bruhat order. Then
$$\partial_{w/u}(fg)=\sum_{w\succeq v\succeq u}u^{-1}v(\partial_{w/v}f)
\partial_{v/u}g \eqno (4.3)
$$
(Generalized Leibnitz' rule).

iii) Let $w=vt$, where $l(w)=l(v)+1$, and $t=t_{ij}$ is the transposition 
that interchanges $i$ and $j$ and fixes all other elements of $[1,n]$. 
Then 
$$\partial_{w/v}=\partial_{ij}, \eqno (4.4)
$$
where $\partial_{ij}:=\partial_{x_ix_j}$.

iv) Let $w_0$ be the longest element of $S_n$. Then
$$w_0v\partial_{w_0/v}=\partial_{w_0v}. \eqno (4.5)
$$

v) Let $u,v,w\in S_n$, $w\succeq u$, and $l(w)=l(u)+l(v)$. Then
$$\partial_{w/u}\s_v=c^w_{uv}, \eqno (4.6)
$$
where $c^w_{uv}$ are the structural constants for the Schubert 
polynomials $\s_w$, $w\in S_n$; in other words,
$$\s_u\s_v\equiv\sum_{w\in S_n}c^w_{uv}\s_w~({\rm mod}~I_n),
$$
where $I_n$ is the ideal generated by the elementary symmetric functions 
$e_1(x_1,\ldots ,x_n),\ldots ,$ $e_n(x_1,\ldots ,x_n)$.}
\vskip 0.2cm

{\it Proof.} We refer the reader to [M1], p.30, for proofs of statements 
i)--iii). 

iv). To prove the identity (4.5), we will use the formula 
(4.2) and the following result due to I.~Macdonald ([M1], (5.7)):
$$\partial_{w_0}(fg)=\sum_{w\in S_n}\epsilon (w)\partial_w(w_0f)
\partial_{ww_0}(g), \eqno (4.7)
$$
where for each permutation $w\in S_n$, $\epsilon (w)=(-1)^{l(w)}$ is the 
sign of $w$.

Using the generalized Leibnitz formula (4.2), we can write the LHS(4.7) 
as follows:
$$\partial_{w_0}(fg)=\sum_{w_0\succeq v}v(\partial_{w_0/v}f)
\partial_v(g). \eqno (4.8)
$$
Comparing the RHS of (4.7) and that of (4.8), we see that
$$v(\partial_{w_0/v}f)=\epsilon (vw_0)\partial_{vw_0}(w_0f).
$$

To finish the proof of equality (4.5), it remains to apply the following 
formula ([M1], (2.12)):
$$\partial_{w_0ww_0}=\epsilon (w)w_0\partial_ww_0.
$$

v). We consider formula (4.6) as a starting point for applications of the skew
divided difference operators to the problem of finding a combinatorial 
formula for the structural constants $c^w_{uv}$ (Littlewood--Richardson
problem for Schubert polynomials). Having in mind some applications of 
(4.6) (see Sections~7 and 8), we reproduce below the proof of (4.6) 
given by I.~Macdonald, [M1], p.112. 
It follows from Proposition~3 $i)$ that
$$c^w_{uv}=\partial_w(\s_u\s_v)=\sum_{w\succeq v_1}v_1
\left(\partial_{w/v_1}\s_u\right)\partial_{v_1}\s_v.
$$
In the latter sum the only nonzero term appears when $v_1=v$. Hence,
$$c^w_{uv}=\partial_w(\s_u\s_v)=v\partial_{w/v}(\s_u)=\partial_{w/v}(\s_u),
$$
since deg~$\partial_{w/v}(\s_u)=0$.
\qed

It is well--known (and follows, for example, from Proposition~2, i) and 
ii)) that for each $v,w\in S_n$
$$\partial_v\left(\s_w\right)\in{\bf N}[x_1,\ldots ,x_{n-1}].
$$
More generally, we make the following conjecture.
\vskip 0.3cm

{\bf Conjecture 1.} {\it For any $u,v,w\in S_n$,
$$\partial_{w/u}{\s_v}\in{\bf N}[X_n],
$$
i.e. $\partial_{w,u}(\s_v)$ is a polynomial in $x_1,\ldots ,x_n$ with 
nonnegative integer coefficients.}
\vskip 0.2cm

{\bf Examples.} Take $w=s_2s_1s_3s_2s_1\in S_4$, $v=s_2s_1\in S_4$, and 
${\bf a}=(2,1,3,2,1)\in R(w)$. There are three possibilities to choose 
${\bf b}$ such that ${\bf b}\subset{\bf a}$, ${\bf b}\in R(v)$, namely, ${\bf 
b}=(2,1,\cdot ,\cdot ,\cdot )$,\break ${\bf b}=(2,\cdot ,\cdot ,\cdot ,1)$ 
and ${\bf b}=(\cdot ,\cdot ,\cdot ,2,1)$. Hence,
$$\partial_{w/v}=s_1s_2s_2s_1\partial_3\partial_2\partial_1+
s_1s_2s_2\partial_1\partial_3\partial_2s_1+
s_1s_2\partial_2\partial_1\partial_3s_2s_1=
\partial_3\partial_2\partial_1-\partial_1\partial_3\partial_{13}-
\partial_{13}\partial_2\partial_{14}.
$$
Using this expression for the divided difference operator 
$\partial_{w/v}$, one can find

a) $\partial_{w/v}(x^3_1x^2_2)=x_1^2+x_1x_4+x_4^2\equiv x_2x_3~({\rm 
mod}~I_4)$. Thus, 
$$\partial_{w/v}(x_1^3x_2^2)\equiv\s_{23}-\s_{13}+\s_{21}~({\rm 
mod}~I_4).
$$ 
We used the following formulae for Schubert polynomials:
$$\s_{23}=x_1x_2+x_1x_3+x_2x_3,~~\s_{13}=x_1^2+x_1x_2+x_1x_3,~~\s_{12}=x_1^2.
$$

b) $\partial_{w/v}(x_1^3x_2^2x_3)\equiv x_2^2x_3~({\rm mod}~I_4)$, and
$$\partial_{w/v}(x_1^3x_2^2x_3)\equiv\s_{121}+\s_{232}-\s_{123}-\s_{213}
-\s_{312}~({\rm mod}~I_4).
$$

c) $\partial_{13}(x_1^3x_2x_3)=x_1x_2x_3(x_1+x_3)\equiv 
-x_1x_2^2x_3~({\rm mod}~I_4)$.
Let us note that $x_1^3x_2x_3=\s_{12321}$.

These examples show that in general

$\bullet$ the "intersection" numbers $\langle\partial_{w/v}(\s_u),
\s_{\tau}\rangle_0$ may  have negative values;

$\bullet$ coefficients $c_{\alpha}$ in decomposition
$\partial_{w/v}(\s_u)\equiv
\ds\sum_{\alpha\subset\delta_n}c_{\alpha}x^{\alpha}~({\rm mod}~I_n)
$  
may take  negative values.
\vskip 0.5cm

{\bf \S 5. Analog of skew divided differences in the Bracket algebra.}
\vskip 0.3cm

In this Section for each $v,w\in S_n$ we construct the element $[w/v]$ 
in the Bracket algebra 
${\cal E}_n^0$ which is an analog of the skew divided difference operators 
$\partial_{w/v}$. The Bracket algebra ${\cal E}_n^0$ was introduced 
in [FK2]. By 
definition, the Bracket algebra ${\cal E}_n^0$ (of type $A_{n-1}$) is the 
quadratic algebra (say, over ${\bf Z}$) with generators $[ij]$, $1\le 
i<j\le n$, which satisfy the following relations
$$\eno{
(i)~~~~&[ij]^2=0,~~{\rm for}~~i<j;\cr
(ii)~~~&[ij][jk]=[jk][ik]+[ik][ij]\cr
&[jk][ij]=[ik][jk]+[ij][ik],~~{\rm for}~~i<j<k;\cr
(iii)~~&[ij][kl]=[kl][ij]\cr
&{\rm whenever}~~\{ i,j\}\cap\{ k,l\} =\emptyset,~~i<j~~{\rm and}~~k<l.}
$$
For further details, see [FK2], [K1].

Now, let $v,w\in S_n$, and $w\succeq v$ with respect to the Bruhat order 
on $S_n$. Let ${\bf a}\in R(w)$ be a reduced decomposition for $w$. We 
define the element $[w/v]$ in the Bracket algebra ${\cal E}_n^0$ to be
$$\eno{
&[w/v]=v^{-1}\sum_{{\bf b}\subset{\bf a},~~{\bf b}\in R(v)}
\phi ({\bf a},b),~~~{\rm where}\cr
&\phi({\bf a},{\bf b})=\prod_i\phi_i({\bf a},{\bf b}),~~{\rm and}\cr
&\phi_i=\cases{s_{a_i}, & $a_i\in{\bf b}$,\cr [a_i,a_i+1], & 
$a_i\not\in{\bf b}$.}}
$$
\vskip 0.3cm

{\bf Remark.} Let $w,v\in S_n$, and $w\succeq v$. One can show that the 
element $[w/v]\in{\cal E}_n^0$ is independent of the reduced decomposition 
${\bf a}\in R(w)$.
\vskip 0.3cm

{\bf Conjecture 2.} {\it The element $[w/v]\in{\cal E}_n^0$ can be written 
as a linear combination of monomials in the Bracket algebra with nonnegative 
integer coefficients.}
\vskip 0.3cm

{\bf Example.} Take $w=s_2s_1s_3s_2s_1\in S_4$, $v=s_2s_1\in S_4$. Then
$$\eno{
[w/v]&=[34][23][12]-[12][34][13]-[13][23][14]\cr
&=[34][12][13]+[34][13][23]-[12][34][13]-[13][23][14]\cr
&=[13][14][23]+[14][34][23]-[13][23][14]\cr
&=[14][34][23].}
$$
\vskip 0.5cm

{\bf \S 6. Skew Schubert polynomials.}
\vskip 0.3cm

{\bf Definition 5.} {\it Let $v,w\in S_n$, and $w\succeq v$ with respect 
to the Bruhat order. The skew 
Schubert polynomial $\s_{w/v}$ is defined to be
$$\s_{w/v}=\partial_{v^{-1}w_0/w^{-1}w_0}(x^{\delta_n}). \eqno (6.1)
$$
}
\vskip 0.3cm

{\bf Examples.} a) Let $w=s_1s_2s_3s_1\in S_4$, and $v=s_1\in S_4$. 
Then $v^{-1}w_0=s_2s_1s_3s_2s_1$, $w^{-1}w_0=s_2s_1$, and
$$\s_{w/v}=\partial_{21321/21}(x_1^3x_2^2x_3)=
(x_1^2+x_1x_4+x_4^2)x_2\equiv\s_{121}+\s_{232}-\s_{124}
-\s_{213}-\s_{321}~({\rm mod}~I_4).
$$

b) Take $w=s_3s_2\in S_4$ and $v=s_3\in S_4$. Then 
$v^{-1}w_0=s_1s_2s_1s_3s_2$, $w^{-1}w_0=s_1s_2s_3s_2$, and 
$\partial_{v^{-1}w_0/w^{-1}w_0}=\partial_{13}$. Thus
$$\s_{w/v}=x_1^2x_2x_3(x_2+x_3)\equiv -x_1^3x_2x_3~({\rm mod}~I_4).
$$

It is clear that if $w,v\in S_n$, and $w\succeq v$, then $\s_{w/v}$ is a 
homogeneous polynomial of degree $\ds\pmatrix{n\cr 2}-l(w)+l(v)$ with 
integer coefficients. It would be a corollary of Conjecture~1 that skew 
Schubert polynomials have in fact positive integer coefficients.

\vskip 0.3cm
{\bf Proposition 4.} {\it

i) Let $v\in S_n$, and $w_0\in S_n$ be the longest element. 
Then 
$$\s_{w_0/v}=\s_{v}. \eqno (6.2)
$$

ii) Let $w\in S_n$, then}
$\s_{w/1}=w_0ww_0\s_{ww_0}.$
\vskip 0.2cm

Proof of (6.2) follows from (4.5) and (6.1).
\qed

It is an interesting task to find the Monk formula for skew Schubert 
polynomials, in other words, to describe the decomposition of the product 
$(x_1+\cdots +x_r)\s_{w/v}$, $w,v\in S_n$, $1\le r\le n-1$, in terms of 
Schubert polynomials.

\vskip 0.5cm
{\bf \S 7. Proof of Conjecture 1 for divided difference operators 
$\partial_{ij}$.}
\vskip 0.3cm

First of all we recall the 
definition of the nilCoxeter algebra $NC_n$ and construction of the 
Schubert expression $\s^{(n+1)}\in{\bf N}[x_1,\ldots ,x_n][NC_n]$. The 
study of action of divided difference operators
$\partial_{ij}$, $1\le i<j\le n$, 
on the Schubert expression $\s^{(n+1)}$ is the main step of 
our proof of Conjecture~1 for the skew divided difference operators
corresponding to the edges in the Bruhat order on the symmetric group 
$S_{n+1}$. In exposition we follow to [FS] and [FK].
\vskip 0.3cm

{\bf Definition 6.} {\it The nilCoxeter algebra $NC_n$ is the algebra 
(say, over ${\bf Z}$) with generators $e_i$, $1\le i\le n$, which satisfy 
the following relations}
$$\eno{
(i)~~~~&e_i^2=0,~~{\rm for}~~1\le i\le n,\cr
(ii)~~~&e_ie_j=e_je_i,~~{\rm for}~~1\le i,j\le n,~|i-j|>1,\cr
(iii)~~&e_ie_je_i=e_je_ie_j,~~{\rm for}~~1\le i,j\le n,~|i-j|=1.}
$$
For each $w\in S_{n+1}$ let us define $e_w\in NC_n$ to be 
$e_w=e_{a_1}\cdots e_{a_p}$, where $(a_1,\ldots ,a_p)$ is any reduced 
decomposition for $w$. The elements $e_w$, $w\in S_{n+1}$, are 
well--defined and form a ${\bf Z}$--basis in the nilCoxeter algebra $NC_n$.

Now we are going to define the Schubert expression $\s^{(n+1)}$ which 
is a noncommutative generating function for Schubert polynomials. Namely,
$$\s^{(n+1)}=\sum_{w\in S_{n+1}}\s_we_w\in{\bf N}[x_1,\ldots ,x_n][NC_n].
$$
The basic property of the Schubert expression $\s^{(n+1)}$ is that it 
admits the following factorization ([FS]):
$$\s^{(n+1)}=A_1(x_1)\cdots A_n(x_n), \eqno (7.1)
$$
where $A_i(x)=\ds\prod_{j=n}^i(1+xe_j)=(1+xe_n)(1+xe_{n-1})\cdots 
(1+xe_i)$.

Now we are ready to formulate and prove the main result of this Section, namely,
the following positivity theorem:
\vskip 0.3cm

{\bf Theorem 1.} {\it Let $1\le i<j\le n+1$, $w\in S_{n+1}$. Then}
$$\partial_{ij}\s_w\in{\bf N}[x_1,\ldots ,x_{n+1}].
$$
\vskip 0.2cm

{\it Proof.} It is enough to prove Theorem~1 for the transposition 
$(i,j)=(1,n)$. 
Thus, we are going to prove that $\partial_{1n}\s_w\in{\bf N}[x_1,\ldots 
,x_n]$. For this goal, let us consider the Schubert expression 
$\s^{(n+1)}=A_1(x_1)A_2(x_2)\cdots A_n(x_n)$, see (7.1). We are 
going to prove that
$$\partial_{1n}\s^{(n+1)}=\sum_{w\in S_{n+1}}\alpha_w(x)e_w,
$$
where $\alpha_w(x)\in{\bf N}[x_1,\ldots ,x_n]$ for all $w\in S_{n+1}$. 
Using the Leibnitz rule (2.2), we can write
$$\eno{
\partial_{1n}\s^{(n+1)}&=\partial_{1n}(A_1(x_1)A_2(x_2)\cdots A_n(x_n))\cr
&=\partial_{1n}(A_1(x_1))A_2(x_2)\cdots A_n(x_n)+A_1(x_n)A_2(x_2)\cdots 
A_{n-1}(x_{n-1})\partial_{1n}(A_n(x_n)).}
$$
First of all, $\partial_{1n}A_n(x_n)=\ds{1+x_ne_n-1-x_1e_n\over x_1-x_n}
=-e_n$. The next observation is
$$\partial_{1n}A_1(x_1)={A_1(x_n)-1\over x_n}+f(x_1,x_n),
$$
where $f(x_1,x_n)\in{\bf N}[x_1,x_n][NC_n]$. Indeed, if 
$A_1(x)=\ds\sum_{k=1}^{n-1}c_kx^k$, where $c_k\in NC_n$, $c_0=1$, then
$$\partial_{1n}A_1(x_1)=\sum_{k=1}^{n-1}c_k{x_1^k-x_n^k\over x_1-x_n}=
\sum_{k=1}^{n-1}c_kx_n^{k-1}+f(x_1,x_n),
$$
and $f(x_1,x_n)\in{\bf N}[x_1,x_n][NC_n]$, as it was claimed. Hence,
$$\eno{
x_n\partial_{1n}\s^{(n+1)}&=(A_1(x_n)-1)A_2(x_2)\cdots A_{n-1}(x_{n-1})
(1+x_ne_n)\cr
&-A_1(x_n)A_2(x_2)\cdots A_{n-1}(x_{n-1})e_nx_n+
F(x_1,\ldots ,x_n)\cr
&=A_1(x_n)A_2(x_2)\cdots A_{n-1}(x_{n-1})-A_2(x_2)A_3(x_3)\cdots A_n(x_n)
+F(x_1,\ldots ,x_n),}
$$
where $F(x_1,\ldots ,x_n)\in{\bf N}[x_1,\ldots ,x_n][NC_n]$. Thus, it is 
enough to prove that the difference
$$A_1(x_n)A_2(x_2)\cdots A_{n-1}(x_{n-1})-A_2(x_2)A_3(x_3)\cdots A_n(x_n)
$$
belongs to the set ${\bf N}[x_1,\ldots ,x_n][NC_n]$. We will use the 
following result (see [FS], [FK1]):
$$A_i(x)A_i(y)=A_i(y)A_i(x),~~~1\le i\le n.
$$
Thus, using a simple observation that $A_i(x)=A_{i+1}(x)(1+xe_i)$, we 
have
$$\eno{
&A_1(x_n)A_2(x_2)\cdots A_{n-1}(x_{n-1})=A_2(x_n)(1+x_ne_1)A_2(x_2)\cdots 
A_{n-1}(x_{n-1})\cr
&=A_2(x_n)A_2(x_2)\cdots A_{n-1}(x_{n-1})+x_nA_2(x_n)e_1A_2(x_2)\cdots 
A_{n-1}(x_{n-1})\cr
&=A_2(x_2)A_2(x_n)A_3(x_3)\cdots 
A_{n-1}(x_{n-1})+x_nA_2(x_n)e_1A_2(x_2)\cdots A_{n-1}(x_{n-1})\cr
&=A_2(x_2)A_3(x_n)(1+x_ne_2)A_3(x_3)\cdots A_{n-1}(x_{n-1})+
x_nA_2(x_n)e_1A_2(x_2)\cdots A_{n-1}(x_{n-1})\cr
&=\cdots =A_2(x_2)A_3(x_3)\cdots A_n(x_n)+x_n\sum_{i=1}^{n-1}\prod_{j=2}^i
A_j(x_j)A_{i+1}(x_n)e_i\prod_{j=i+1}^{n-1}A_j(x_j). & (7.2)}
$$
Let us denote  the sum over $i$ in (5.3) by $F(x_1,\ldots ,x_n)$. It is 
clear that 
$$F(x_1,\ldots x_n)\in{\bf N}[x_1,\ldots x_n][NC_n].
$$ 
Thus the difference 
$$A_1(x_n)A_2(x_2)\cdots A_{n-1}(x_{n-1})-A_2(x_2)A_3(x_3)\cdots A_n(x_n)
=F(x_1,\ldots ,x_n)
$$
also belongs to the set ${\bf N}[x_1,\ldots ,x_n][NC_n]$.

\qed
\vskip 0.5cm

{\bf \S 8. Generating function for Schubert polynomials structural 
constants $c_{uv}^w$.}
\vskip 0.3cm

Let $w,v\in S_n$, $l(w)-l(v)\le 1$, and $w\succeq v$ with respect to 
the Bruhat order. For $1\le i\le n$ and $1\le s\le n-1$ we define the 
element $e_i^{(s)}(w/v)$ of the nilCoxeter algebra 
$NC_n$ using the following rule
$$e_i^{(s)}(w/v)=\cases{0, &if $w=vt_{(a,b)}$, and simultaneously 
$a\ne s$ and $b\ne s$,\cr e_{n-i}, &if $w=vt_{(s,b)}$, and $s<b$, \cr
-e_{n-i}, & if $w=vt_{(b,s)}$, and $b<s$,\cr
1, & if $w=v$.}
$$
\vskip 0.3cm

{\bf Theorem 2.} {\it Let $u,w\in S_n$. Then
$$\sum_{v\in S_n}c_{uv}^we_v=\sum_{\{ v_i^{(s)}\}^{n-1}_{s=1}}
\prod_{s=1}^{n-1}\prod_{i=1}^{n-s}e_i^{(s)}\left( 
v_{i-1}^{(s)}/v_i^{(s)}\right) , \eqno (8.1)
$$
summed over all sequences of permutations $\{ v_0^{(s)}\succeq v_1^{(s)}
\succeq\cdots\succeq v_{n-s}^{(s)}\}^{n-1}_{s=1}$ such that
$$v_0^{(1)}=w, v_0^{(s+1)}=v_{n-s}^{(s)},~1\le s\le n-2,~v_1^{(n-1)}=u.
$$
In the product in the RHS(8.1) the factors are multiplied 
left--to--right, according to the increase of $s$.}
\vskip 0.2cm

{\it Proof.} We start with rewriting of the LHS(8.1), namely, $\ds\sum_{v\in S_n}
c_{uv}^we_v=\eta\left(\partial_{w/u}\s^{(n)}\right)$, where $\s^{(n)}$ 
denote the Schubert expression. Indeed,
$$\eta\left(\partial_{w/v}\s^{(n)}\right)=\sum_{v\in S_n}\eta
\left(\partial_{w/v}(\s_v)\right) e_v=\sum_{v\in S_n}c_{uv}^we_v.
$$
The next step is to compute $\eta\left(\partial_{w/v}\s\right)$ using the 
generalized Leibnitz rule (4.3) and the following Lemma:
\vskip 0.3cm

{\bf Lemma 1.} {\it Let $w,u\in S_n$, and $f_1,\ldots ,f_N\in P_n$. Then}
$$u\partial_{w/u}(f_1\cdots f_N)=\sum_{w=v_0\succeq v_1\succeq\cdots
\succeq v_{N-1}\succeq v_N=u}~\prod_{i=1}^Nv_i\left(\partial_{v_{i-1}/v_i}
(f)\right) .
$$
\vskip 0.2cm

We apply Lemma~1 to the Schubert expression
$$\s^{(n)}=A_1(x_1)\cdots A_{n-1}(x_{n-1})=\prod_{i=1}^{n-1}\prod_{k=n-1}^i
(1+x_ie_k). \eqno (8.2)
$$
In the RHS(8.2) in the product the factors are multiplied left--to--right 
according to the increase of $i$. 
As a result, we obtain
$$\eta\left(\partial_{w/u}\s^{(n)}\right)=\sum_{\{ v_0^{(s)}\succeq v_1^{(s)}
\succeq\cdots\succeq v_{n-s}^{(s)}\}_{s=1}^{n-1}}\prod_{s=1}^{n-1}
\prod_{i=1}^{n-s}\eta\left(\partial_{v_{i-1}^{(s)}/v_i^{(s)}}
(1+x_se_{n-i})\right) ,
$$
summed over all sequences of permutations $\{ v_0^{(s)}\succeq
v_1^{(s)}\succeq\cdots\succeq v_{n-s}^{(s)}\}_{s=1}^{n-1}$ such that 
$v_0^{(1)}=w$, $v_0^{(s+1)}=v_{n-s}^{(s)}$, $1\le s\le n-2$, 
$v_1^{(n-1)}=u$.

It is clear that we can assume $l(v_{i-1}^{(s)})-l_(v_i^{(s)})\le 1$ for 
all $i,s$, and under these conditions, we have
$$\eta\left(\partial_{v_{i-1}^{(s)}/v_i^{(s)}}(1+x_se_{n-i})\right)=
e_i^{(s)}(v_{i-1}^{(s)}/v_i^{(s)}).
$$
\qed

\vskip 0.5cm
{\bf \S 9. Open problems.}
\vskip 0.3cm

Below we formulate a few problems related to the content of this paper.

1) (Main problem) Let $w,v\in S_n$ and $w\succeq v$ with respect to the 
Bruhat order on the symmetric group $S_n$. To prove that polynomials 
$\partial_{w/v}\left(\s_u\right)$ have nonnegative coefficients for each 
$u\in S_n$.

2) (Generalized Littlewood--Richardson problem for Schubert polynomials)

Let $u,v,w\in S_n$ and
$$\partial_{w/v}\left(\s_u\right)=\sum_{\alpha}c_{\alpha}x^{\alpha}.
$$
To find a combinatorial description of coefficients 
$c_{\alpha}:=c_{\alpha}(u,v,w)$.

\vskip 0.3cm
{\bf Remark.} If $l(w)=l(u)+l(v)$ and $w\succeq v$, then $\partial_{w/v}
\left(\s_u\right)=c^w_{uv}$, see [M1], p.112, or the present paper, 
Proposition~3, $v)$.
\vskip 0.2cm

3) (Skew key polynomials) Let $\alpha =(\alpha_1,\ldots ,\alpha_n)$ be
a composition, $\lambda (\alpha )$ be the unique partition in the orbit 
$S_n\cdot\alpha$, and $w(\alpha )\in S_n$ be the shortest permutation 
such that $w(\alpha )\cdot\alpha =\lambda (\alpha )$. Let $v\in S_n$ 
be such that $w(\alpha )\succeq v$ with respect to the Bruhat order. Using
in Definition~4 the isobaric divided difference operators 
$\pi_i:=\partial_ix_i$, $1\le i\le n-1$  (see, e.g., [M1], p.28) instead
of operators $\partial_i$ one can define for each pair $w\succeq v$
 the skew isobaric divided difference operator $\pi_{w/v}:P_n\to P_n$.
We define the skew key polynomial $k_{\alpha /v}$ to be
$$k_{\alpha /v}=\pi_{w(\alpha )/v}\left( x^{\ld(\alpha )}\right) ,
$$
where $x^{\beta}=x_1^{\beta_1}\cdots x_n^{\beta_n}$ for any composition
$\beta =(\beta_1,\beta_2,\ldots ,\beta_n)$. It is natural to ask whether 
or not the skew key polynomials have nonnegative coefficients?

4) To find a geometrical interpretation of the skew divided difference 
operators, the polynomials $\partial_{w/v}\left(\s_u\right)$, and the 
skew key polynomials.

5) Does there exist a stable analog of the skew Schubert polynomials?

\vskip 0.5cm
{\bf References.}
\vskip 0.5cm

\item{[BS]} Bergeron N. and Sottile F., {\it Skew Schubert functions and 
the Pieri formula for flag manifolds,} Preprint alg-geom/9709034;

\item{[FK1]} Fomin S. and Kirillov A.N., {\it The Yang--Baxter equation, 
symmetric functions and Schubert polynomials,} Discrete Mathematics, 
1996, v.53, p.123-143;

\item{[FK2]} Fomin S. and Kirillov A.N., {\it Quadratic algebras, Dunkl 
elements, and Schubert calculus,} Preprint AMSPPS\#199703--005--001, 
1997, 34p.; to appear in {\it Progress in Geometry,} ed. J.-L.~Brylinski 
and R.~Brylinski, 1997;

\item{[FS]} Fomin S. and Stanley R., {\it Schubert polynomials and 
nilCoxeter algebra,} Adv. in Math., 1994, v.103, p.196-207;

\item{[K1]} Kirillov A.N., {\it On some quadratic algebras,} 
Preprint CRM-2478, 1997, 28p; 

q-alg/9705003.

\item{[K2]} Kirillov A.N., {\it Lectures on Schubert polynomials,} 
      in preparation;
      
\item{[LS]} Lascoux A. and Sch\"utzenberger M.-P., {\it Polyn\^omes de Schubert,}
C.R. Acad. Sci. Paris, Serie~I, 1982, t.294, p.447-450;

\item{[M1]} Macdonald I.G., {\it Notes on Schubert polynomials}, Publ.
 LACIM, 1991, Univ. du Quebec \`a Mont\'eal;

\item{[M2]} Macdonald I.G., {\it Symmetric functions and Hall 
polynomials,} Second ed., Oxford Univ. Press, New York/London, 1995;

\end